# ON THE ORIGIN OF THE ARABIC NUMERALS

A. BOUCENNA

Département de Physique, Faculté des Sciences,
Université Ferhat Abbas 19000 Sétif, Algérie

aboucenna@wissal.dz

**Abstract**

From the pagination of an Algerian Arabic manuscript of the beginning of the 19th century, we rediscover the original shape that the Arabic numerals had before passing in Europe and underwent the transformation that gave the modern Arabic numerals. This original shape, whose use disappeared completely, proves that these numerals have their origin in the Arabic letters. Contrary to what some hypotheses pretend, particularly those that present them as drifting of Indian characters, the 10 Arabic numerals that we use are, nothing else, 10 Arabic letters more or less modified and taken in the "Abjadi" order. The hypothesis of the Indian origin of the Arabic numerals is revealed a mistake denied by the shape of the Arabic numerals and by the logic of the right to left representation of the numbers and the algorithm of the elementary operations. The Arabic numerals that simplified the writing of the numbers and the algorithms of the elementary operations are believed to be born in the Maghreb (North Africa). From Béjaïa (Bougie) they passed in to Europe to give, after evolution, the modern Arabic numerals : 0, 1, 2, 3, 4, 5, 6, 7, 8, 9. they also migrate to the Middle East (the Mashrek) to give, after transformations in shape and adding two Hebrew letters, the Arabic numerals, "Mashrekis", that are used currently in Middle East : ٠ ١ ٢ ٣ ٤ ٥ ٦ ٧ ٨ ٩ .

# SUR L'ORIGINE DES CHIFFRES ARABES

A. BOUCENNA
Département de Physique, Faculté des Sciences,
Université Ferhat Abbas 19000 Sétif, Algérie
aboucenna@wissal.dz

**Résumé**

A travers la pagination d'un manuscrit arabe algérien du début du 19$^{ème}$ siècle nous redécouvrons la forme originale que les chiffres arabes avaient avant de passer en Europe et subir les transformations qui ont donné les chiffres arabes modernes. Cette forme originale, dont l'utilisation a complètement disparu, montre que ces chiffres ont pour origine les lettres arabes. Contrairement à ce que prétendent certaines hypothèses, particulièrement celles qui les présentent comme dérivant de caractères indiens, les 10 chiffres arabes que nous utilisons ne sont en fait que 10 lettres arabes plus ou moins modifiées et données dans l'ordre "Abjadi". L'hypothèse de l'origine indienne des chiffres arabes se révèle une erreur démentie par la forme des chiffres arabes et par la logique de la sociologie droite-gauche de la représentation des nombres et des algorithmes des opérations élémentaires. Les chiffres arabes qui ont beaucoup simplifié l'écriture des nombres et les algorithmes des opérations élémentaires sont nés au Maghreb (Afrique du nord). De Béjaïa (Bougie) ils sont passés en Europe pour donner après évolution les chiffres arabes modernes : 0, 1, 2, 3, 4, 5, 6, 7, 8, 9. Ils sont aussi passés au Moyen Orient (le Mashrek) pour donner après transformations, remise en forme et en ajoutant deux lettres hébraïques, les chiffres arabes "Mashreki" utilisés actuellement au Moyen Orient : ٠ ١ ٢ ٣ ٤ ٥ ٦ ٧ ٨ ٩ .

**Mots clés :** manuscrit arabe algérien, chiffre arabe original, chiffre arabe moderne, chiffre "Mashreki", lettres arabes, lettres hébraïques, ordre "Abjadi", logique de la sociologie droite-gauche.

## I. Introduction

L'origine des chiffres arabes modernes : 0, 1, 2, 3, 4, 5, 6, 7, 8, 9 qui se sont imposés dans le monde entier, a fait l'objet de plusieurs hypothèses. Une de ces hypothèses les présente comme dérivant de caractères indiens [1], ce qui a amené certains à les qualifier de chiffres indo-arabes. Une autre hypothèse fait une relation directe entre le nombre d'angles plausibles dans la forme géométrique du chiffre est sa valeur [2]. Le chiffre arabe quantifierait le nombre d'angles de la forme géométrique qui le représente. Ce raisonnement semble avoir comme point de départ la forme moderne des chiffres arabes. L'origine doit être recherchée dans la forme originale qu'avaient les chiffres arabes avant de passer en Europe. L'évolution importante que ces chiffres ont subie par la suite n'a rien à voir avec leur origine. A travers la pagination d'un manuscrit arabe algérien du début du 19$^{ème}$ siècle nous redécouvrons la forme originale du chiffre arabe dont l'utilisation a complètement disparu. Cette forme montre de façon évidente que les chiffres arabes ont pour origine des lettres arabes et prouve que l'hypothèse de l'origine indienne des chiffres arabes est démunie de tout fondement. Les chiffres arabes que nous connaissons actuellement et qui ont beaucoup simplifié la représentation des nombres et les algorithmes des opérations élémentaires sont nés au Maghreb. De Béjaïa, ils sont passés en Europe. Le chiffre arabe ou chiffre "Ghoubar" a été toujours lié à la notion de "Ghoubari" qui signifie 'calcul'. Cette dénomination montre l'importance de l'impulsion donnée par les chiffres arabes à la simplification des algorithmes de calcul.

**II. Forme originale des chiffres arabes (les chiffres arabes du manuscrit)**

**1. Le Manuscrit**

Le manuscrit, de format 17.5 cm par 12.5 cm et dans lequel les chiffres arabes sous leur forme originale ont été utilisés, est une copie du célèbre livre "Kitab khalil bni Ishak El Maliki". Il a 179 feuilles numérotées de 1 à 179, et lui manque les feuilles 4 et 23. La feuille 77 a été numérotée, certainement par erreur, deux fois. Elle porte le numéro 77 au recto et le numéro 78 au verso. Ce livre traite du "fik'h" islamique selon le rite malikite, très répandu au Maghreb. C'est une référence en la matière. Des copies sont disponibles dans toutes les Zouïas, et tous les anciens Instituts. Notre manuscrit est une copie réalisée par Mohamed ben Sidi khélifa ben Ali ben Salem ben Sifa, pour le Compte de la Zaouïa de "Sidi Khélifa" de la région de Sétif en Algérie. La date de cette copie est donnée à la fin du manuscrit, par la l'expression traditionnelle : 'Copie achevée à Midi ("Ezzaouel") du Mercredi, huit du mois de "Safar" de l'Année 1225 de l'ère hégire', correspondant approximativement à l'année 1810 après J.C. (Figure 1a). La date 1810 n'est pas très ancienne. Mais il ne faut pas considérer cette datte en valeur absolue mais il faut la considérer en valeur relative. Dans cette région du monde, en 1800, les gens ignoraient tout de ce qui se passait en Europe. Ils ignoraient tout de la révolution industrielle, des philosophes, des mathématiciens, des savants européens. Pour eux la référence c'est toujours les savants grecs, c'est aussi Ibnou Sina, c'est Ibnou Roushd, c'est El Farabi, c'est Ibn Khaldoun, etc…, mais pas Descartes, ni Avogadro, ni Pascal,… 1810 signifie en fait à peu près 1200 ou 1400, voire 1000.

**Figure 1a**
Partie de la dernière page du Manuscrit
Kitab Khalil bni Ishak

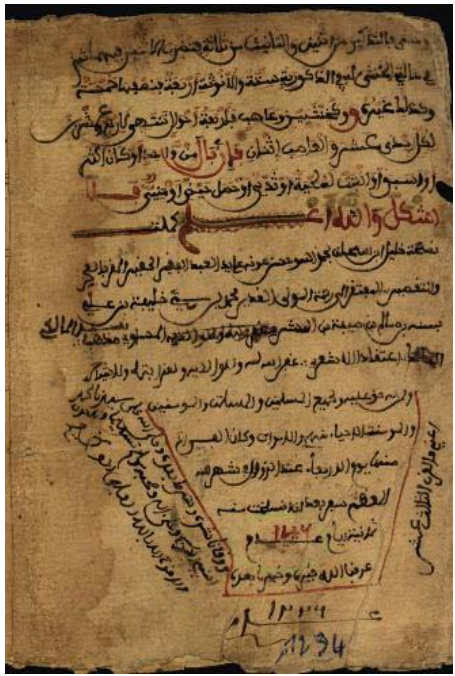

L'année du manuscrit est indiquée en chiffres arabes originaux en dessus de l'expression عــــــام (entourée par une ellipse) et qui signifie année. Un essai de transcription, venu certainement plus tard, de cette datte en chiffre "Mashreki" n'a pas été fructueux. Puisque le chiffre 5 a été pris pour un 6 "Mashreki". Un autre essai de transcription mais en chiffres arabes modernes en bas de la page a été lui aussi soldé par un échec puisque le chiffre 5 du manuscrit a été interprété comme étant le chiffre 4. Ceci montre la confusion totale à laquelle été confrontés les Maghrébins de la deuxième moitié du 19[ème] siècle, qui venait de perdre leurs chiffres arabes originaux en faveur des chiffres "Mashreki" et des chiffres arabes modernes.

Dans ce manuscrit l'année est justement exprimée en chiffres arabes mais sous leur forme originale (Figure 1b). Cette présentation des dates des copies de manuscrits en chiffres arabes originaux est retrouvée dans des copies d'autres manuscrits de la même époque (Figures 2 et 3). Les Maghrébins ont certainement continué à utiliser les chiffres arabes sous leur forme originale, indépendamment de l'évolution subie par ces chiffres en Europe jusqu'au début du 19$^{ème}$ siècle. Il faut signaler que juste après cette époque (1810) arrivèrent des évènements très douloureux liés à l'occupation coloniale. Les repères changent. L'identité maghrébine résiste et se tourne vers le "Mashrek" pour trouver un allié. Dès lors les Maghrébins oublièrent leurs chiffres, les chiffres arabes originaux, et adoptèrent les chiffres "Mashreki" du Moyen Orient (le "Mashrek") (Figure 4). D'ailleurs on constate un moment d'incertitude et de transition, où apparemment on ne sait plus quels sont les chiffres qu'il faut utiliser (figure 5).

**Figure 1b**
Partie de la dernière page du Manuscrit Kitab Khalil bni Ishak
daté en chiffres arabes originaux de 1225 de l'ère Hégire

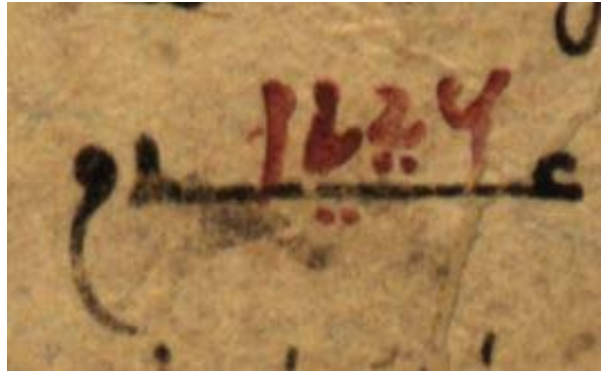

**Figure 2**
Partie de la dernière page d'un Manuscrit
daté en chiffres arabes originaux de 1149 de l'ère Hégire

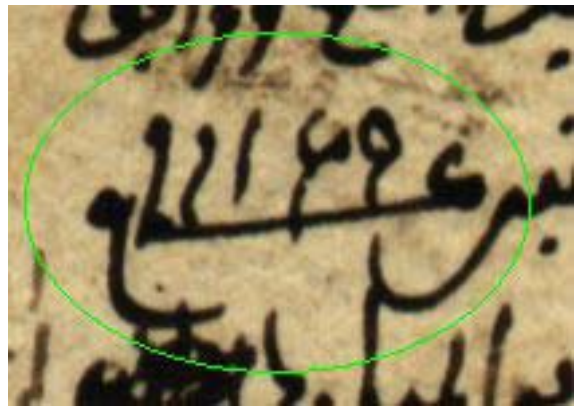

**Figure 3**
Partie de la dernière page d'un autre Manuscrit
daté en chiffres arabes originaux de 1194 de l'ère Hégire

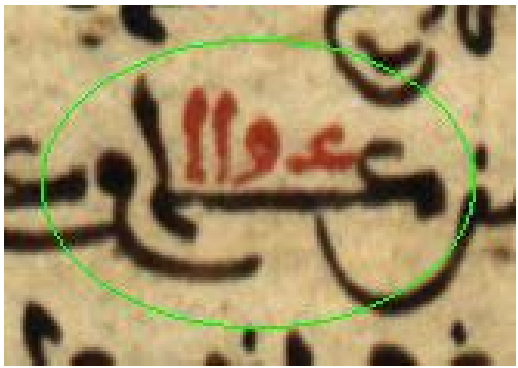

**Figure 4**
Partie de la dernière page d'un autre Manuscrit
daté en chiffres "Mashreki" de 1311 de l'ère Hégire

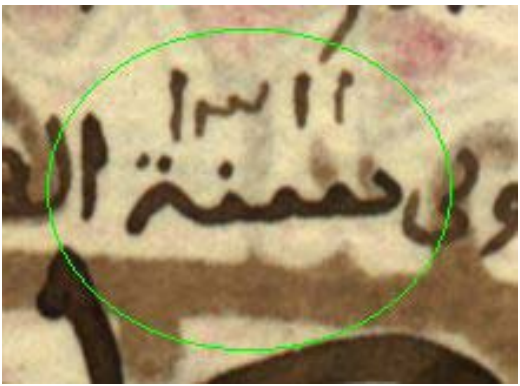

**Figure 5**
Partie de la dernière page d'un autre Manuscrit

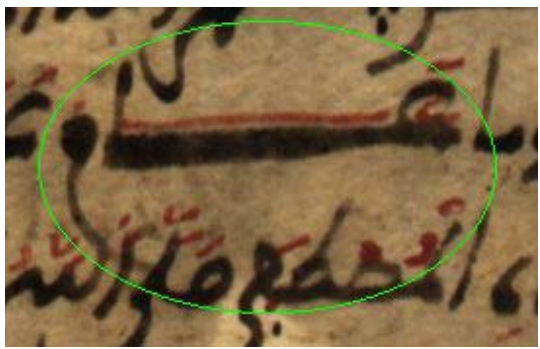

L'année de cette copie n'est pas mentionnée, l'emplacement traditionnellement réservé est laissé vide. Alors que le jour, le mois voire l'heure sont mentionnés. C'est vraiment l'embarras.

## 2. Pagination mixte et identification des chiffres du manuscrit

Dans les manuscrits arabes, il est rare de trouver l'ordre des pages indiqué par une suite de nombres, comme est le cas de la pagination moderne des ouvrages. C'est la pagination traditionnelle qui était largement utilisée. Elle consiste à inscrire, au pied gauche de la page paire de la feuille, le premier mot ou les premiers deux mots, s'il y a ambiguïté, de la première ligne de la page impaire de la feuille suivante (Figure 6).

**Figure 6**
Pagination mixte

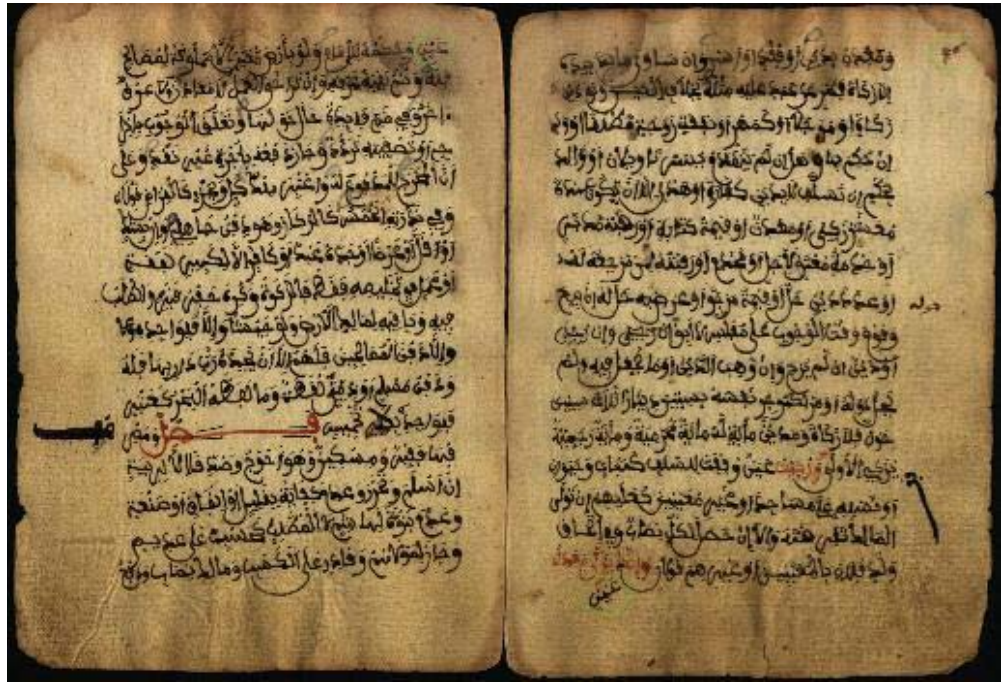

En haut à droite de la page droite, est indiqué le numéro de la feuille 36. En bas à gauche de la même feuille est indiqué le mot عين (entouré par une ellipse) qui est le premier mot de la première ligne de la page suivante, la page gauche.

Notre manuscrit est un spécimen rare, puisque en plus de la pagination traditionnelle, la numérotation des feuilles a été employée en utilisant justement les chiffres arabes originaux. L'ordre des pages est ainsi indiqué, dans chaque feuille, de deux manières différentes. L'intérêt de cette **pagination mixte** est double. D'abord nous avons la forme initiale, originale que les chiffres arabes eurent, avant de passer en Europe et avant d'être modifiés pour donner la forme actuelle, moderne des chiffres arabes. Ensuite la pagination mixte et le nombre important de pages du manuscrit éliminent tout risque d'erreur dans l'interprétation et l'identification des symboles. Nous avons 179 feuilles toutes numérotées. Chaque chiffre est donc cité dans plus de 18 suites ordonnées.

## 3. Forme originale des chiffres arabes

Dans la suite de ce texte, nous qualifierons de forme **moderne** la forme actuelle des chiffres arabes et de forme **originale** la forme qu'eurent les chiffres arabes avant de passer en Europe et de subir les transformations ayant abouti à la forme actuelle, moderne des chiffres arabes. Les figures 7, 8, 9, 10, 11, 12, 13, 14, 15 et 16 montrent les chiffres utilisés pour la numérotation des feuilles de notre manuscrit. C'est cette forme originale des chiffres arabes qui était en vigueur dans les pays du Maghreb jusqu'au début du 19^(ème) siècle. L'identification des symboles utilisés dans le manuscrit pour représenter les chiffres, montrés sur les figures 7 à 16, et reportés dans la table 1, est déduite de la pagination mixte du manuscrit.

**Figure 7**
De gauche à droite
Les chiffres arabes 1, 0 et 0 du manuscrit

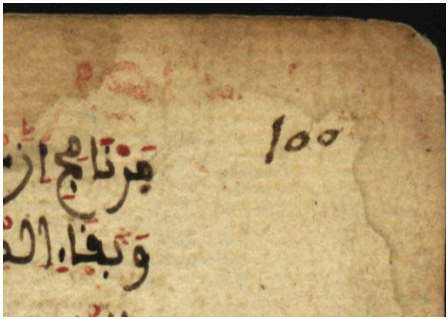

**Figure 8**
de gauche à droite
Les chiffres arabes 2 et 1 du manuscrit

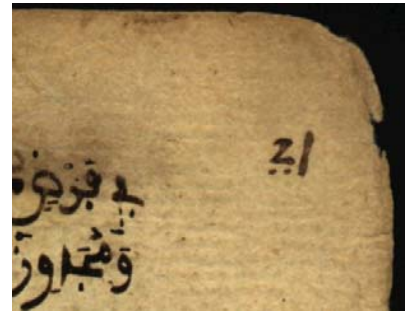

**Figure 9**
De gauche à droite
Les chiffres arabes 2 et 2 du manuscrit

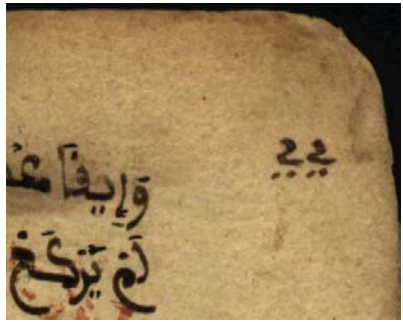

**Figure 10**
De gauche à droite
Les chiffres arabes 3 et 3 du manuscrit

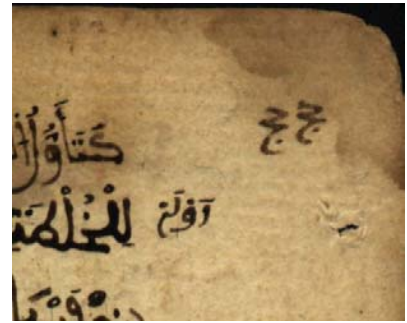

**Figure 11**
De gauche à droite
Les chiffres arabes 2 et 4 du manuscrit

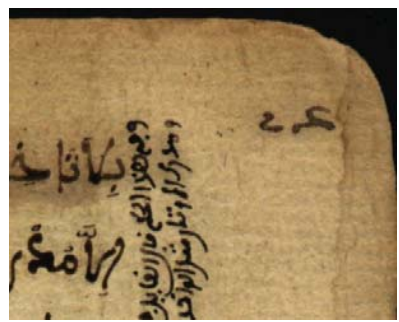

**Figure 12**
De gauche à droite
Les chiffres arabes 3 et 4 du manuscrit

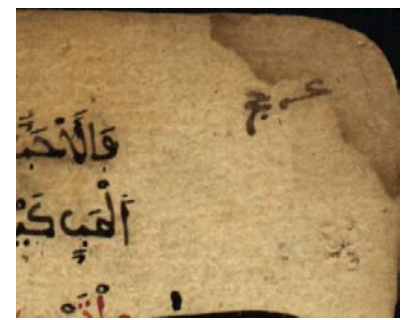

**Figure 13**
De gauche à droite
Les chiffres arabes 3 et 5 du manuscrit

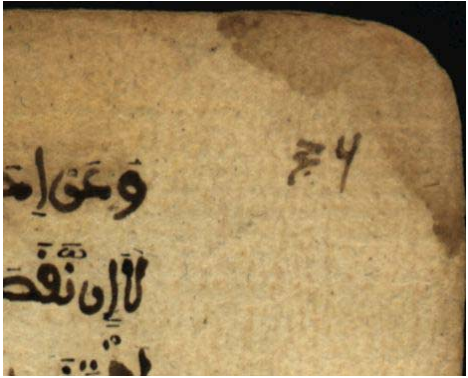

**Figure 14**
De gauche à droite
Les chiffres arabes 3 et 7 du manuscrit

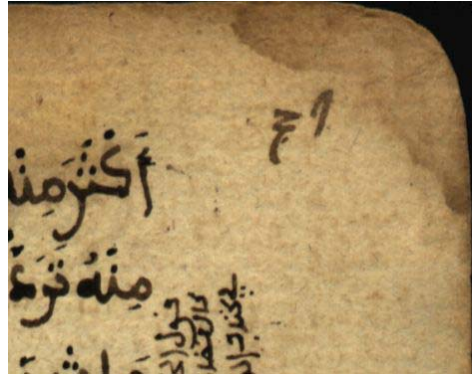

**Figure 15**
De gauche à droite
Les chiffres arabes 6 et 8 du manuscrit

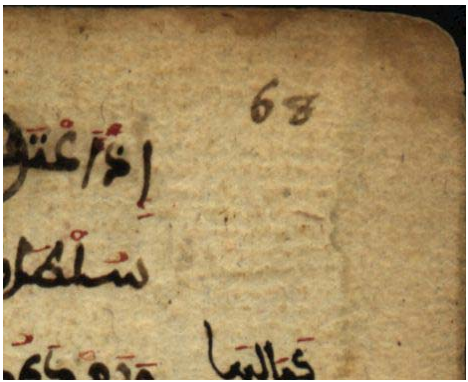

**Figure 16**
De gauche à droite
Les chiffres arabes 6 et 9 du manuscrit

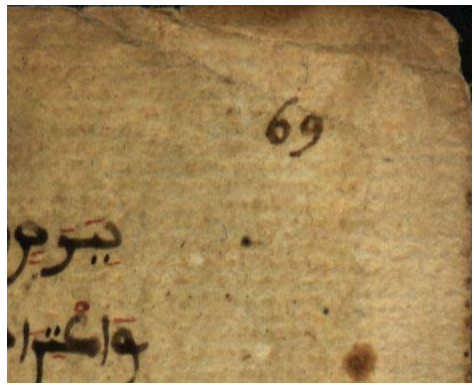

### 4. Correspondances

Hormis certaines divergences apparentes minimes, on constate :
- d'abord une nette ressemblance entre les formes des chiffres du manuscrit et les formes des chiffres arabes modernes, la correspondance est immédiate,
- ensuite une ressemblance évidente entre la forme des chiffres du manuscrit et la forme "Maghrebi" de certaines lettres arabes de l'ordre "Abjadi".

La forme "Maghrebi" des lettres arabes était très répandue au Maghreb. En langage moderne ce serait la fonte "Maghrebi" des caractères arabes. La table 1 donne la correspondance entre les chiffres du manuscrit, les chiffres arabes modernes et certaines lettres arabes de l'ordre "Abjadi". Lors du passage des chiffres arabes originaux aux chiffres arabes modernes on constate la permutation des chiffres 4 et 5 et la disparition des points des chiffres 2 et 3 qui rappelaient la relation avec les lettres arabes.

**Table 1**
Correspondance entre les Chiffres Arabes du Manuscrit
et les Lettres Arabes de forme "Maghrébi"

| Chiffres Arabes | Chiffres du Manuscrit | Lettres Arabes | Sons | Formes Maghrébi des Lettres Arabes | | | |
|---|---|---|---|---|---|---|---|
| 1 | ١ | ا | Alif | ) | | ى | |
| 2 | | ب | Baa | ب | | | |
| 2 a) | ج | ي | Yaa | ى | ب | | |
| 3 | ج | ج | Jim | ج | ج | | |
| 4 | د | د | Del | د | | | |
| 5 | ٧ | ه | Haa | ه | ح | ه | ه |
| 6 | 6 | و | Waw | و | | | |
| 7 | ر | ز | Zin | ز | ن | | |
| 8 | 8 | ح | H'aa | ح | | | |
| 9 | و | ط | T'aa | ط | | | |
| 0 | O | ص | Sad | ص | ص | | |

a) La lettre ي ne correspond pas à la valeur 2 de la suite "Abjadi", elle correspond à la valeur 10 (voir texte).

## III. Origines des Chiffres Arabes : Fondements théoriques et arguments

La table 1 montre une ressemblance évidente entre la forme des chiffres arabes originaux et la forme "Maghrebi" de certaines lettres arabes. Cette ressemblance montre que les chiffres arabes ont pour origine les lettres arabes. Maintenant, quelle était la logique qui a motivée le choix de la lettre arabe qui devait représenter tel ou tel chiffre ? La réponse est contenue dans la notion de l'ordre "Abjadi" des lettres arabes.

### 1. Ordre "Abjadi" et valeurs "Abjadi" des lettres

#### 1.1. La suite "Abjadi" des nombres

La **suite "Abjadi" des nombres** est la suite finie des 28 nombres naturels suivants : {1, 2, 3, 4, 5, 6, 7, 8, 9, 10, 20, 30, 40, 50, 60, 70, 80, 90, 100, 200, 300, 400, 500, 600, 700, 800, 900, 1000}. Les nombres appartenant à cette suite sont des nombres "Abjadi"

#### 1.2. Les lettres arabes

L'alphabet d'une langue est énuméré selon un ordre conventionnel. L'alphabet arabe se présente selon deux ordres différents : l'ordre dit "Abjadi" et l'ordre dit "Arabe". L'ordre dit "Arabe" est celui utilisé dans l'apprentissage de la langue arabe, c'est celui utilisé dans la majorité des dictionnaires arabes.

**Table 2**
Ordre "Abjadi" et Valeurs "Abjadi des Lettres Arabes

| Ordre "Abjadi" | Lettre Arabe | Son | Valeur "Abjadi" | Ordre "Abjadi" | Lettre Arabe | Son | Valeur "Abjadi" |
|---|---|---|---|---|---|---|---|
| 1 | ا | Alif | 1 | 15 | س | Sin | 60 |
| 2 | ب | Baa | 2 | 16 | ع | A'in | 70 |
| 3 | ج | Jim | 3 | 17 | ف | Faa | 80 |
| 4 | د | Del | 4 | 18 | ص | Sad | 90 |
| 5 | ه | Haa | 5 | 19 | ق | K'af | 100 |
| 6 | و | Waw | 6 | 20 | ر | Raa | 200 |
| 7 | ز | Zin | 7 | 21 | ش | Shin | 300 |
| 8 | ح | H'aa | 8 | 22 | ت | Taa | 400 |
| 9 | ط | T'aa | 9 | 23 | ث | Thaa | 500 |
| 10 | ي | Yaa | 10 | 24 | خ | Kh'aa | 600 |
| 11 | ك | Kef | 20 | 25 | ذ | Dhel | 700 |
| 12 | ل | Lem | 30 | 26 | ض | Dzad | 800 |
| 13 | م | Mim | 40 | 27 | ظ | Dzaa | 900 |
| 14 | ن | Noun | 50 | 28 | غ | Ghin | 1000 |

L'ordre dit "Abjadi" est basé sur le "calcul Abjadi" dit "Hissab el-joummel". Dans cet **ordre "Abjadi"** les 28 lettres arabes sont quantifiées. A chaque lettre est associée une **valeur "Abjadi"** de la **suite "Abjadi"** des nombres. Dans l'ordre "Abjadi" les lettres sont rangées par ordre croissant de leurs valeurs "Abjadi". La table 2 donnes les valeurs "Abjadi" des lettres arabes rangées dans **l'ordre "Abjadi"**. De nos jours et pour des motivations informatiques, l'idée de la quantification des caractères en établissant une correspondance entre des lettres et des nombres est devenue essentielle. En code ASCII, par exemple, un caractère est représenté par un nombre compris entre 0 et 255.

## 1.3. Les lettres hébraïques

L'alphabet arabe n'est pas le seul alphabet se présentant dans l'ordre "Abjadi". D'autres alphabets des anciennes langues du Moyen Orient se présentent certainement sous cet ordre. L'alphabet hébraïque est rangé de la même façon que l'alphabet arabe. Dans l'**ordre "Abjadi"** des lettres hébraïques, les 22 lettres sont quantifiées. A chaque lettre est associée une **valeur "Abjadi"** de la **suite "Abjadi"** des nombres. Dans l'ordre "Abjadi" des lettres hébraïques les lettres sont rangées par ordre croissant de leurs valeurs "Abjadi". La table 3 donne les valeurs "Abjadi" des lettres hébraïques rangées dans **l'ordre "Abjadi"**.

**Table 3**
Ordre "Abjadi" et Valeurs "Abjadi des Lettres Hébraïques

| Ordre "Abjadi" | Lettre Hébraïque | Son | Valeur "Abjadi" | Ordre "Abjadi" | Lettre Hébraïque | Son | Valeur "Abjadi" |
|---|---|---|---|---|---|---|---|
| 1 | א | Aleph | 1 | 15 | ס | Samek | 60 |
| 2 | ב | Beth | 2 | 16 | ע | Ayin | 70 |
| 3 | ג | Gimel | 3 | 17 | ף פ | Fe | 80 |
| 4 | ד | Daleth | 4 | 18 | ץ צ | Tsahde | 90 |
| 5 | ה | He | 5 | 19 | ק | Q'oph | 100 |
| 6 | ו | Vav | 6 | 20 | ר | Regh | 200 |
| 7 | ז | Zayin | 7 | 21 | ש ש | Sin Shin | 300 |
| 8 | ח | Cheth | 8 | 22 | ת | Tav | 400 |
| 9 | ט | Teth | 9 | | | | |
| 10 | י | Yodh | 10 | | | | |
| 11 | ך כ | Kaph | 20 | | | | |
| 12 | ל | Lamedh | 30 | | | | |
| 13 | ם מ | Mem | 40 | | | | |
| 14 | ן נ | Nun | 50 | | | | |

## 1.4. Utilisation des valeurs "Abjadi"

### 1.4.1. Les chiffres

Grâce à leurs valeurs "Abjadi" les lettres arabes ou hébraïques servaient aussi de chiffres pour représenter les nombres (Figure 17). Pour représenter le nombre 1245, par exemple, les arabes utilisaient le mot "همرغ" composé des lettres arabes : هـ qui a la valeur "Abjadi" 5, م qui la valeur "Abjadi" 40, ر qui a la valeur "Abjadi" 200 et غ qui a la valeur "Abjadi" 1000. Les valeurs "Abjadi" des lettres arabes sont traditionnellement utilisées en poésie arabe pour indiquer les dates des poèmes. Par exemple si un poème date de 1245, on trouve cette date intégrée de façon très intelligente à la fin du poème à travers le mot sans signification linguistique "همرغ". De très habiles commerçants utilisent à nos jours, en Algérie, les valeurs "Abjadi" des lettres arabes pour indiquer les prix des articles commercialisés. Les prix sont représentés par des mots que seuls ces commerçants savent déchiffrer en utilisant justement les valeurs "Abjadi" des lettres arabes qui composent ces mots.

**Figure 17**
Ecriture des nombres en en utilisant les lettres et leurs valeurs "Abjadi"

| | |
|---|---|
| همرغ | 1245 |
| ر | 200 |
| طصظغ | 1945 |

Grâce à leur valeur "Abjadi" les lettres arabes et hébraïques servaient de chiffres pour représenter les nombres.

### 1.4.2. Digitalisation d'expressions ou "Hisseb el-joummel"

Une phrase entière peut être évaluée en calculant la somme des valeurs "Abjadi" des lettres composant les mots de cette phrase. Par exemple, la valeur de l'expression "احمد زينب" est : 
$1 + 8 + 40 + 4 + 7 + 10 + 50 + 2 = 122$.
Cette forme de digitalisation d'expressions utilisant les valeurs "Abjadi" des lettres est très utilisée en astrologie. Les valeurs des expressions sont soumises à des algorithmes de calculs bien déterminés et le résultat final est comparé à des tables astrologiques.

### 2. Stratégie du choix des symboles représentant les chiffres

Tout nombre compris entre 1 et 1999 peut être représenté par un ensemble de lettres arabes (un mot) en utilisant les valeurs "Abjadi" des lettres. De même tout nombres compris entre 1 et 499 peut être représenté par un ensemble de lettres hébraïques (un mot). Au-delà, la représentation des nombres se complique. La même difficulté est rencontrée avec l'utilisation des symboles (chiffres) romains. Pour simplifier les choses, une logique radicalement différente fut adoptée. Elle consistait, probablement, à :
- réduire le nombre de symboles à 10, au lieu de 28 lettres pour représenter les nombres.
- considérer que la longueur des mots représentant les nombres inférieurs 1000 est exactement trois symboles,
- adopter un sens, de droite à gauche, dans l'écriture et la lecture des mots qui doivent représenter les nombres de sorte que le symbole le plus à droite représente les unités, le symbole suivant représente les dizaines et le symbole suivant représente les centaines.

Cette façon de procéder met en évidence le rôle de la notion de rangs : le rang des unités, le rang des dizaines et le rang des centaines et marque, par la même occasion, et surtout, la naissance de la philosophie du chiffre dans sa conception actuelle. La notion de rang étant incluse dans l'ancienne représentation des chiffres par des lettres en considérant leurs valeurs

"Abjadi". Nous avons (Tables 2 et 3) des lettres qui représentent les unités (1, 2, …, 9) des lettres qui représentent les dizaines (10, 20, …, 90) et des lettres qui représentent les centaines (100, 200, …, 900). Dans le cadre de cette approche, l'écriture des nombres serait basée sur l'emploi de nouveaux symboles, les chiffres.

### 2.1. Les lettre-chiffres

Les grecs utilisaient leurs lettres pour représenter les nombres. Les romains utilisaient les lettres latines I, V, X, C, D, M comme chiffres. Les arabes et les hébreux utilisaient leurs lettres pour représentaient les nombres, etc … Avant l'invention des chiffres arabes, avec la philosophie du chiffre dans sa conception actuelle, l'utilisation des lettres de l'alphabet pour représenter les nombres était une pratique courante. C'était l'ère des lettre-chiffres. Les symboles représentant les chiffres étaient de facto des lettres de l'alphabet.

### 2.2. Base 16 et ses chiffres

Aujourd'hui, pour choisir les chiffres permettant de représenter les nombres dans une base donnée, la base 16 par exemple, qui a 16 chiffres, on utilise les dix chiffres arabes de la base 10 et on complète la liste par les lettres de l'alphabet de la civilisation dominante, i.e. les lettres de l'alphabet latin. Les chiffres de la base 16 sont : 0, 1, 2, 3, 4, 5, 6, 7, 8, 9, A, B, C, D, E, F. Personne n'est allé chercher des chiffres illisibles d'une pièce archéologique provenant d'une civilisation utilisant les bases 15 ou 60. Quoique, maintenant, avec les moyens de communication, et les progrès de la science, on pourrait facilement trouver des chiffres appartenant à des civilisations utilisant plus que 10 chiffres pour représenter les nombres dans la base 16.

### 2.3. Base 10 et chiffres arabes

### 2.3.1. Les chiffres 1, 2, 3, 4, 5, 6, 7, 8 et 9

En considérant la stratégie du choix des symboles représentant les chiffres qui était et qui est toujours en vigueur, les chiffres arabes, que nous utilisons actuellement, ont certainement une relation avec les lettres de l'alphabet de la civilisation qui les a inventés. Si les chiffres arabes moderne sont inventés au Maghreb à l'époque florissante de la civilisation arabo-musulmane, les symboles représentant les chiffres arabes 1, 2, 3, 4, 5, 6, 7, 8 et 9 de la base 10 doivent tout naturellement découler des formes des 9 premières lettres arabes sous leur forme "Maghrebi" de l'ordre "Abjadi". Dans la table 4 sont donnée les dix première lettres arabes de l'ordre "Abjadi", leurs valeurs "Abjadi", leurs formes "Maghrebi", les symboles qui seraient les chiffres arabes directement issus des lettres arabes de l'ordre "Abjadi", Les chiffres arabes originaux, les chiffres arabes originaux passés en Europe et les chiffres arabes modernes. Hormis certaines détails, les ressemblances évidentes entre chiffres arabes originaux, ceux du manuscrit, et les 10 première lettres arabes sous leur forme "Maghrebi" données dans l'ordre "Abjadi" sont trop importante pour être de simples coïncidences fortuites.

### 2.3.2. Le chiffre 0 ou "Sifr"

Si le choix des symboles représentant les chiffres arabes 1, 2, …, 9 peut se justifier en partie par les valeurs "Abjadi" des lettres arabes choisies pour représenter ces chiffres, le symbole représentant le chiffre 0 ne peut pas s'expliquer de la sorte, la valeur zéro n'étant pas une valeur "Abjadi".

**Figure 18**
Le sens du mot arabe "Sifr" et de ses dérivées

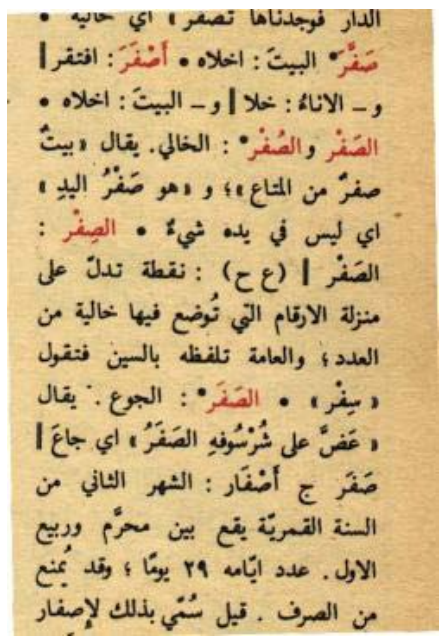

Les babyloniens connaissaient un 0, ils l'utilisaient dans la représentation des nombres, mais il n'était pas utilisé pour simplifier les algorithmes de calcul. La découverte du 0, dans sa conception actuelle, celle que nous connaissons, qui a simplifié la représentation des nombres et les algorithmes des opérations élémentaires a été certainement un fait très marquant. Ce 0 a eu toujours les mêmes missions depuis sa découverte et est ainsi transmis en Europe. Ce 0 est né au Maghreb. Comme tout nouveau né il fallait lui choisir un nom et un symbole. Le nom choisi est le mot "Sifr". Est-t-il un mot arabe ? La figure 18 montre une partie d'une page d'un dictionnaire de la langue arabe qui donne le sens du mot "sifr" et de ses dérivés. On apprend que le mot "Sifr" n'est pas un mot étranger à la langue arabe. Il n'est pas emprunté à une autre langue pour décrire un état nouveau méconnu par la langue arabe, comme est le cas par exemple, du mot "Falsafa" qui désigne "philosophie" et qui est emprunté à la langue grecque pour décrire un état nouveau. Le sens du mot "sifr" et de ses dérivés décrit toujours un état vide auquel on ne s'attendait pas. Le vide décrit par le mot "sifr" est dans tous les cas une situation anormale. Le nom "Sifr" donné au 0 exprime exactement le rôle que doit jouer le 0 dans la représentation des nombres par des chiffres. Le 0 doit remplacer un chiffre qui a fait défaut dans un rang donné.

Le symbole choisi pour représenter le "Sifr" ou "Zéro" est le symbole 0 ? Ce symbole est-t-il une lettre arabe ? Si oui, Laquelle ? Revenons à notre logique du choix des symboles représentant les chiffres. La valeur 0 n'est pas une valeur "Abjadi" il serait difficile de trouver une lettre arabe pour représenter le 0 en ne considérant que les valeurs "Abjadi" des lettres. La stratégie de choix des symboles doit être modifiée pour ne pas tenir compte des valeurs "Abjadi". On pourrait alors choisir la dixième lettre arabe pour représenter le 0, c'est-à-dire la lettre arabe ي (Yaa) correspondant au nombre 10 de la suite "Abjadi". Ce n'est pas le cas. A mon avis, c'est la forme initiale la lettre arabe ص (Sad) qui a servi à fabriquer le chiffre arabe original 0. Le choix de la lettre arabe ص n'a pas été motivé par les valeurs "Abjadi" des lettres arabes. La lettre ص (Sad) est tout simplement la première lettre du mot صفر "Sifr" qui désigne le chiffre 0 en arabe. Le 0 ou "Sifr" dans sa conception actuelle, par le sens de son nom "Sifr" et par son symbole 0 est une invention de la civilisation arabo-musulmane. D'ailleurs les mots 'Chiffre' et 'Zéro' passés par la suite dans les autres langues ont été dérivés du mot "Sifr". Au Maghreb, où les chiffres arabes, sous leur forme originale, jouaient pleinement leur rôle dans le calcul, le "Ghoubari", tous les chiffres s'appelaient probablement "Asfar". Par la suite au Mashrek (Moyen Orient) fut utilisé le mot "Rakm" pour désigner le chiffre.

### IV. L'ère du chiffre

#### 1. Logique des transformations initiales

Comme la valeur 0 n'est pas une valeur "Abjadi", le chiffre arabe 0 a été le premier chiffre formé comme tel, indépendamment de l'ordre "Abjadi" des lettres arabes et de leurs valeurs "Abjadi". L'idée de fabriquer les autres chiffres 1, 2, …, 9 indépendamment des lettres est certainement, venue juste après pour consacrer définitivement le divorce entre les chiffres et les lettres avec leurs valeurs "Abjadi". En effet l'ordre "Abjadi" aurait voulu que le symbole représentant le chiffre 5 soit la forme isolée de la lettre ه (Haa), cinquième lettre de l'ordre "Abjadi" de l'Alphabet arabe correspondant au nombre 5 de la suite "Abjadi". Or cette espèce de cercle est déjà utilisée pour désigner le 0. Il a fallu trouver un autre symbole pour désigner le chiffre 5. De cette priorité on déduit que le chiffre arabe 0 a été le premier chiffre découvert. L'invention des chiffres arabes a été réalisée en un temps très court. Une fois l'idée du chiffre découverte, Il a fallut peut être 1 à 2 jours pour former tous les chiffres arabes y compris le 0. Ce qui explique la liberté prise par le Maître du chiffre au moment de la mise en forme de tous les chiffres. Il n'a pas pris la lettre arabe ب (Baa), comme l'aurait voulu la suite "Abjadi" pour désigner le chiffre 2, mais il a choisi la forme finale "Maghrebi" de la lettre arabe ي (Yaa) avec ses deux points. Mieux encore, il s'est même permis d'omettre les deux points. L'arabe s'écrit avec des lettres attachées. Pour pouvoir représenter les chiffres, le Maître du chiffre a fait subir aux lettres arabes les transformations nécessaires visant à écarter le caractère attaché des lettres arabes. Ainsi il n'a pas pris directement la lettre و (Waw) comme le voulait l'ordre "Abjadi" pour représenter le chiffre 6, mais il a décidé de prendre la forme renversée de la lettre و (Waw). Il n'a pas pris directement la lettre ط (T'aa) comme le voulait l'ordre "Abjadi" pour représenter le chiffre 9, mais il a pris la forme renversée de la lettre ط (T'aa). Notre Maître a eu la même attitude lors de la mise en forme de tous les chiffres. Cette liberté montre bien la volonté affichée de se séparer complètement de la vieille représentation des nombres qui était intimement liée aux valeurs "Abjadi" des lettres arabes. C'est une nouvelle ère qui commence et notre Maître était apparemment bien conscient de l'importance de sa réalisation. Il a agi librement, avec un esprit nouveau, mais sans oublier complètement l'ordre "Abjadi" des lettres arabes. C'est cette démarche qui fait que l'origine de la forme des chiffres arabes est bien la forme des lettres arabes de la suite "Abjadi" mais que leur esprit était complètement différent de l'esprit des valeurs "Abjadi".

#### 2. Mise en forme des chiffres

#### Le Chiffre 0

Le chiffre arabe 0 a été en fait le premier vrai chiffre découvert. Son symbole est la forme initiale de la lettre arabe ص (Sad). Le choix de la lettre arabe ص n'a pas été motivé par les valeurs "Abjadi" des lettres arabes, comme a été le cas pour les autres chiffres arabes. La lettre ص (Sad) est tout simplement la première lettre du mot صفر "Sifr" qui désigne le chiffre 0 en arabe.

**Le chiffre 1**

Le chiffre 1 du manuscrit (figure 8) est exactement la lettre ا (Alif), première lettre de l'ordre "Abjadi" de l'alphabet arabe correspondant au nombre 1 de la suite "Abjadi".

**Le Chiffre 2**

En observant la forme du chiffre 2 du manuscrit avec ses deux points (figures 7, 8, 9 et 11), on fait rapidement la relation avec la lettre arabe ي (Yaa) correspondant au nombre 10 de la suite "Abjadi". C'est exactement une des formes finale, isolée de la forme "Maghrebi" de la lettre arabe ي. Dans le manuscrit les deux points sont parfois omis.

**Le Chiffre 3**

Le Chiffre 3 du manuscrit (figures 10, 12, 13 et 14) est la forme finale liée, légèrement modifiée de la lettre ج (Jim) de l'ordre "Abjadi" de l'alphabet arabe correspondant au nombre 3 de la suite "Abjadi". On remarque la présence du point du ج (Jim) dans le chiffre 3 du manuscrit. Pour rompre avec la possibilité de l'écrire en attaché, la lettre ج (Jim) a été dédoublée par sa forme initiale, mais sans le point. Cette forme originale, classique du chiffre 3, celle du manuscrit, avec sont point, est un argument de taille en faveur du fait que les chiffres arabes ont pour origine les lettres arabes.

**Le Chiffre 4**

Le chiffre 4 du manuscrit (figures 11 et 12) est exactement la forme Maghrebi de la lettre د (Del), quatrième lettre de l'ordre "Abjadi" de l'alphabet arabe correspondant au nombre 4 de la suite "Abjadi". C'est la forme originale du chiffre 4, celle du manuscrit, qui argumente le plus le fait que les chiffres arabes ont pour origine les lettres arabes.

**Le Chiffre 5**

Le chiffre 5 du manuscrit (figure 13) est une forme modifiée de la forme finale liée de la forme Maghrebi de la lettre ھ (Haa), cinquième lettre de l'ordre "Abjadi" de l'Alphabet arabe correspondant au nombre 5 de la suite "Abjadi".

**Le chiffre 6**

Le Chiffre 6 du manuscrit (figures 15 et 16) est la lettre و (Waw) renversée. C'est la sixième lettre de l'ordre "Abjadi" de l'Alphabet arabe correspondant au nombre 6 de la suite "Abjadi". Cette transformation veut rompre avec la possibilité d'écrire en attaché la lettre و (Waw).

**Le Chiffre 7**

Le chiffre 7, du manuscrit (figure 14) est une forme légèrement modifiée de la forme Maghrebi de la lettre ز (Zin), septième lettre de l'ordre "Abjadi" de l'Alphabet arabe et correspondant au nombre 7 de la suite "Abjadi". On constate que le point a été remplacé par un petit trait collé au corps de la lettre. Cette transformation veut éviter l'écriture en attaché de la lettre ز (Zin).

**Le chiffre 8**

Le chiffre 8, du manuscrit (figure 15) est une forme légèrement modifiée de la forme finale de la forme Maghrebi de la lettre ح (H'aa), huitième lettre de l'ordre "Abjadi" de l'Alphabet arabe et correspondant au nombre 8 de la suite "Abjadi". Pour éviter d'écrire en attaché la lettre ح (H'aa) une barre oblique vient lier la queue de la lettre à son début.

**Le Chiffre 9**

Le chiffre 9 (figure 16) est une forme renversée de la lettre ط (T'aa), neuvième lettre arabe de l'ordre "Abjadi" et correspondant au nombre 9 de la suite "Abjadi". Cette transformation vise à rompre avec la possibilité d'écrire en attaché la lettre ط (T'aa) et de rompre avec l'ancienne représentation des nombres.

**V. Evolution des chiffres**

Les chiffres arabes originaux nés au Maghreb, ont émigré vers le Mashrek pour devenir, après évolution, les chiffres "Mashreki" et vers l'Europe pour devenir après évolution, les chiffres arabes modernes. Notons que les arabes du Mashrek, les "Mashreki", qualifient les chiffres arabes du Maghreb, de chiffres "Franji", des Francs, et les arabes du Maghreb qualifient les chiffres "Mashreki", de chiffres "Hindi". Dans les deux cas les qualificatifs "Franji" et "Hindi" n'ont rien à voir avec l'origine des chiffres arabes du Maghreb et des chiffres arabes du Mashrek. C'est peut être ce qualificatif "Hindi" qui a induit en erreur certains chercheurs allant trouver des origines indiennes aux chiffres arabes. Ces chercheurs ignorent certainement tout des formes originales des chiffres arabes et de l'ordre "Abjadi" des lettres arabes. Dans la sociologie Maghrebine le qualificatif "Hindi" ne signifie pas forcément d'origine indienne, il peut signifier tout simplement : magique, bizarre, surprenant, exotique. Les fruits appelés en France **figues de barbarie** (du Maghreb ancien) s'appellent dans leur propre pays, le Maghreb, **"Hindi"** certainement à cause de leur apparence assez exotique.

**1. Les chiffres arabes modernes**

Hormis la permutation des chiffres 4 et 5, la comparaison des formes modernes et originales des chiffres arabes montre que le chiffre arabe a pratiquement gardé sa forme initiale. Les améliorations qui ont été apportées touchaient surtout le côté esthétique. La permutation des chiffres 4 et 5 montre tout simplement qu'il n'a pas été accordé une grande importance à l'ordre "Abjadi" des lettres arabes, ce qui est compatible avec l'esprit du chiffre tel que le voulait le Maître du chiffre. La table 4 montre l'évolution des chiffres arabes originaux vers les chiffres arabes modernes.
Les formes des chiffres arabes modernes 0 et 1 sont pratiquement les mêmes que les formes originales. Dans la forme moderne du chiffre arabe 2, les deux points ont complètement disparu. Dans la forme moderne du chiffre arabe 3, le point et la queue en relation avec la lettre arabe ج (Jim) ont disparu. L'observation des chiffres arabes modernes 4 et 5, et leur comparaison aux chiffres arabes originaux 4 et 5 (Table 1 et figures 11 et 13) montre que les deux représentations ont été permutées. En effet le chiffre arabe classique 5 ressemble plutôt au chiffre arabe moderne 4, alors que le chiffre arabe classique 5 ressemble plutôt au chiffre arabe moderne 4. Même la forme "imprimerie" du chiffre arabe moderne 4 est proche de la forme finale liée de la lettre ـه (Haa), cinquième lettre de l'ordre "Abjadi" de l'Alphabet arabe. Nous sommes certainement en face d'une simple inversion de l'ordre des symboles qui représentent les chiffres arabes 4 et 5. Ceci peut être dû à un acte volontaire ou, tout simplement, à une erreur. La forme moderne du chiffre arabe 6 est pratiquement la même que le forme classique. Une légère évolution est constatée dans la forme moderne du chiffre arabe moderne 7.

**Table 4**
Evolution des formes des Chiffres Arabes originaux
vers les Chiffres Arabes modernes

| Lettre Arabe | Son | Valeurs Abjadi | Lettre Arabe Maghrebi | Chiffres Arabes 1) | Chiffres Arabes 2) | Chiffres Arabes 3) | Chiffres Arabes Modernes |
|---|---|---|---|---|---|---|---|
| ا | Alif | 1 | ۱ | ۱ | ۱ | ۱ | 1 |
| ب | Baa | 2 | ب | ب | ݘ | ݘ | 2 |
| ج | Jim | 3 | ج | ج | ج | ج | 3 |
| د | Del | 4 | ک | ک | ک | 4 | 4 |
| ه | Haa | 5 | ۵ | ۵ | 4 | ک | 5 |
| و | Waw | 6 | و | و | 6 | 6 | 6 |
| ز | Zin | 7 | ز | ز | 7 | 7 | 7 |
| ح | H'aa | 8 | ح | ح | 8 | 8 | 8 |
| ط | T'aa | 9 | 6 | 6 | 9 | 9 | 9 |
| ي | Yaa | 10 | ي | ي | | | |
| ص | Sad | 90 | ص | ص | O | O | 0 |

1) Première version éventuelle des chiffres arabes, respectant scrupuleusement l'ordre "Abjadi" des lettres arabes.
2) Version modifiée des chiffres arabes, tenant compte de la logique des transformations initiales (voir texte).
3) Version européenne initiale des chiffres arabes. Remarquer la permutation des chiffres 4 et 5.

Dans la forme moderne du chiffre arabe 8 les angles sont arrondis. La forme moderne du chiffre arabe 9 ne diffère pas de façon significative de la forme du chiffre 9 du manuscrit.

**Les chiffres "Mashréki"**

Au Moyen-Orient (le Machrek), où le chiffre arabe est qualifié de chiffre "Frangi", c'est le chiffre "Mashréki", qualifié par les Maghrébins de "Hindi", qui est très utilisé. D'où viennent ces chiffres "Mashréki" ?
L'observation des chiffres "Mashréki" (Table 5), particulièrement les chiffres 1, 4, 5 et 9 permet d'affirmer que les chiffres " Mashréki" ont aussi pour origine les chiffres arabes originaux. La comparaison des formes classiques des chiffres arabes et celles des chiffres " Mashréki" montre des changements qui ont voulu rétablir l'emprise des valeurs "Abjadi" des lettres arabes sur le chiffre. Ce qui s'écarte de la volonté de séparer la représentation des nombres en considérant les valeurs "Abjadi" des lettres arabes.

**Table 5**
Evolution des formes des Chiffres Arabes originaux
vers les Chiffres Arabes "Mashreki"

| Valeurs "Abjadi" | Lettre Arabe | Son | Chiffres Arabes Originaux | Chiffres Arabes Modernes | Chiffres Arabes Mashreki | Lettre Hebraique | Son |
|---|---|---|---|---|---|---|---|
| 1 | ا | Alif | ∫ | 1 | ١ | א | Aleph |
| 2 | ب | Baa | ح | 2 | ٢ | ב | Beth |
| 3 | ج | Jim | ح | 3 | ٣ | ג | Gimel |
| 4 | د | Del | ے | 4 | ٤ | ד | Daleth |
| 5 | ه | Haa | 4 | 5 | ٥ | ה | He |
| 6 | و | Waw | 6 | 6 | ٦ | ו | Vav |
| 7 | ز | Zin | 7 | 7 | ٧ | ז | Zayin |
| 8 | ح | H'aa | 8 | 8 | ٨ | ח | Cheth |
| 9 | ط | T'aa | 9 | 9 | ٩ | ט | Teth |
| 10 | ي | Yaa |  |  | . | י | Yadh |
| 90 | ص | Sad | O | 0 |  |  |  |

La stratégie des transformations consistait à :
  - redessiner pratiquement la forme des chiffres en respectant scrupuleusement l'ordre "Abjadi" des lettres, et en commençant par le chiffre 1. Le chiffre 0 n'est plus le chiffre le plus important,
  - se débarrasser du caractère attaché des lettres arabes,
  - remplacés les points des lettres arabes par une sorte de pied collé au corps de la lettre,
  - éviter de prendre des symboles identiques aux chiffres arabes originaux pour représenter des chiffres différents pour empêcher toute confusion,
  - emprunter des lettres hébraïques à la place des lettres arabes chaque fois qu'une confusion avec les chiffres arabes originaux est possible ou qu'une origine "Abjadi" d'un chiffre n'est pas visible comme est le cas pour le chiffre 6 et pour le chiffre 0.

**Le chiffre 1** "Mashréki" ١ (table 5) est exactement la lettre ا (Alif), première lettre de l'ordre "Abjadi" de l'alphabet arabe correspondant au nombre 1 de la suite "Abjadi".

**Le chiffre 2** "Mashréki" ٢ (table 5) est représenté par une forme modifiée de la lettre ب (Baa), deuxième lettre arabe de l'ordre "Abjadi" des lettres arabes et correspondant au nombre 2 de la suite "Abjadi". Le point est remplacé par une sorte de pied collé au corps de la lettre. Le chiffre 2 n'est plus représenté par la lettre arabe ي (Yaa) correspondant au nombre 10 de la suite "Abjadi".

**Le chiffre 3** "Mashréki" ٣ (table 5) est une forme modifiée par une rotation de π/2 à droite du chiffre 3 original qui est lui même la forme finale liée, légèrement modifiée de la lettre ج (Jim) de l'ordre "Abjadi" de l'alphabet arabe correspondant au nombre 3 de la suite "Abjadi". Le point est remplacé par une sorte de pied collé au corps de la lettre.

**Le chiffre 4** "Mashréki" ٤ (table 5) est une forme légèrement déformée de la forme originale du chiffre arabe 4, qui est lui même exactement la forme Maghrebi de la lettre د (Del), quatrième lettre de l'ordre "Abjadi" de l'Alphabet arabe correspondant au nombre 4 de la suite "Abjadi".

**Le chiffre 5** "Mashréki" ٥ (table 5) a été représenté par la forme isolée de la lettre ه (Haa), cinquième lettre de l'ordre "Abjadi" de l'Alphabet arabe et correspondant au nombre 5 de la suite "Abjadi". C'est une forme très proche de celle retenue par le Maître du chiffre pour représenter le chiffre arabe 0

**Le chiffre 6** "Mashréki" ٦ (table 5) a été représenté par la forme isolée de la lettre ו (Vav), sixième lettre de l'ordre "Abjadi" de l'Alphabet hébraïque et correspondant au nombre 6 de la suite "Abjadi". Le choix de la lettre arabe و (Waw) sixième lettre de l'ordre "Abjadi" de l'Alphabet arabe correspondant au nombre 6 de la suite "Abjadi" est empêchée par la confusion possible qu'aurait pu avoir lieu avec le chiffre arabe original 9.

**Le chiffre 7** "Mashréki" ٧ (table 5) est une forme modifiée de la lettre ز (Zin), septième lettre de l'ordre "Abjadi" de l'Alphabet arabe et correspondant au nombre 7 de la suite "Abjadi". Le point est remplacé par une sorte de pied collé au corps de la lettre.

**Le chiffre 8** "Mashréki" ٨ (table 5) est une forme modifiée par une rotation de π/2 à gauche de la lettre ح (H'aa), huitième lettre de l'ordre "Abjadi" de l'Alphabet arabe et correspondant au nombre 8 de la suite "Abjadi".

**Le chiffre 9** "Mashréki" ٩ (table 5) est une forme renversée de la lettre ط (T'aa), neuvième lettre arabe de l'ordre "Abjadi" et correspondant au nombre 9 de la suite "Abjadi". C'est la même que la forme originale du chiffre arabe 9.

**Le chiffre 0** Le symbole ه choisi pour représenter le chiffre 5 "Mashréki" a une forme très proche de celle retenue par le Maître du chiffre pour représenter le chiffre arabe 0. Il a fallut choisir un autre symbole. Une logique de choix des symboles possible consisterait à prendre la lettre arabe ي (Yaa) correspondant au nombre 10 de la suite "Abjadi" comme symbole du chiffre "Mashréki" 0. Ce choix peut entraîner une confusion avec le chiffre arabe original 2.

On emprunta alors la lettre י (Yodh), dixième lettre de l'ordre "Abjadi" de l'Alphabet hébraïque et correspondant au nombre 10 de la suite "Abjadi". Ce choix montre bien que le premier chiffre "Mashréki" qui était mis en forme était le chiffre 1 et non pas le chiffre 0.

### VI. Le chiffre arabe, une révolution

### 1. Les chiffres calculette

Certes l'homme depuis ses débuts a essayé de quantifier la nature à l'aide de nombres. Pour représenter ces nombres il a utilisé des chiffres. Chaque civilisation a eu ses chiffres. Les chiffres arabes, avec leur zéro, ont été une révolution à plus d'un titre. D'abord ils ont simplifié la représentation des nombres, ils ont simplifié les algorithmes de calcul et surtout ils nous ont dotés d'une petite calculette pour faire des opérations très compliquées. Cette calculette n'est pas un instrument mais une méthode simple qui ne nécessite pas plus d'une feuille et d'un crayon. A quelle civilisation appartenait l'homme qui a été à l'origine de cette grande invention. Par l'origine des formes des chiffres arabes, nous venons de montrer que ces 10 chiffres arabes que nous utilisons ne sont en fait que 10 lettres arabes plus ou moins modifiées données dans un ordre bien déterminé qui est l'ordre "Abjadi". Les chiffres en eux-mêmes n'ont d'intérêt que par leurs incidences sur la représentation des nombres et sur la simplification des algorithmes de calcul. Si nous examinons de plus près cette représentation des nombres et ces algorithmes de calcul, nous découvrons une logique que nous appelons la **logique de la sociologie droite-gauche.**

**2. Logique de la sociologie droite-gauche**

Une personne appartenant à une société ayant **une logique de la sociologie droite-gauche** (Figure 19) écrit et lit un mot ou un nombre de droite à gauche. Pour écrire le nombre : 12457892, il commence par écrire les unités, les dizaine les centaines des unités, ensuite les unités, les dizaines, les centaines des milles etc … en commençant par la droite. Il fera de même pour la lecture en commençant par la droite : Il lira : 2 et 90 et 800 ; 7 et 50 et 400 milles ; 2 et 10 millions. Une personne appartenant à une société ayant **une logique de la sociologie gauche-droite** écrit et lit un mot ou un nombre de gauche à droite. Pour lire le même nombre 12457892, qui est écrit avec une logique de la sociologie droite-gauche, qui n'est pas la sienne, Il commence d'abord par se faire des repères de droite à gauche, c'est-à-dire, il essaye d'abord de s'y retrouver dans une logique de la sociologie droite-gauche, une fois qu'il isole les chiffres composant ce nombre de droite à gauche 12 457 892, il revient à sa logique de la sociologie gauche-droite, et lit : 12 millions  457 milles  892.
Toutes les opérations élémentaires dont les chiffres arabes ont simplifié les algorithmes s'effectuent dans une logique de la sociologie droite-gauche. D'ailleurs cette logique nous est imposée, par les chiffres arabes, chaque fois qu'il s'agit de traiter avec des nombres.
Si c'était une civilisation de la logique de la sociologie gauche-droite qui était à l'origine du chiffre arabe, la civilisation indienne, par exemple, comme le prétendent certains, la manipulation des nombres serait forcément de gauche à droite. Pour écrire le nombre 12457892, par exemple, on aurait écrit de gauche à droite les unités ensuite les dizaine puis les centaines etc … c'est-à-dire, écrire 29875421 au lieu de 12457892 et on aurait économisé à un esprit de la logique de la sociologie gauche-droite, l'acrobatie du changement de la logique de la sociologie chaque fois qu'il doit manipuler des nombres. Je ne vois aucune raison qui oblige un inventeur à changer de logique de la sociologie au moment de la réalisation de son invention. Les chiffres arabes portent la signature da la civilisation qui les a inventés qui est la civilisation arabo-musulmane.
Le fait que les chiffres arabes ont pour origine des formes "Maghrebi" des lettres arabes, que les chiffres arabes ont été très utilisés dans le Ghoubar (calcul) au Maghreb (Algérie) jusqu'au début du 19$^{ème}$ siècle prouve que les chiffres arabes, pur produit de la civilisation arabo-musulmane, sont nés au Maghreb (Afrique du nord). Il sont passés en Europe une première fois avec le père Sylvestre mais ils n'ont pas eu de succès. C'est de Béjaïa (Bougie) qu'ils sont passés une deuxième fois en Europe avec des commerçants, et des mathématiciens comme Leonardo Fibonacci, qui les ont appréciés.
En Europe, les chiffres arabes ont donné, après évolution, les chiffres arabes modernes
0, 1, 2, 3, 4, 5, 6, 7, 8, 9. Ils  sont passés aussi au Moyen Orient (le Mashrek) pour donner après transformations, remise en forme et en ajoutant deux lettres hébraïques les  chiffres arabes "Mashreki" utilisés actuellement au Moyen Orient ٠ ١ ٢ ٣ ٤ ٥ ٦ ٧ ٨ ٩ .

**Figure 19**
Logique de la Sociologie

---

**Lecture du nombre : 12457892 ?**

**Personne appartenant à la Logique de la Sociologie Droite-Gauche :**

**Lecture en un temps :**

2 et 90 et 800 et 7 et 50 et 400 mille et 2 et 10 millions

**Personne appartenant à la Logique de la Sociologie Gauche-Droite :**

**Lecture en deux temps :**

**Temps 1 : passage en à Logique de la Sociologie Droite-Gauche :**

12 457 892

**Temps 2 : retour à la Logique de la Sociologie Gauche-Droite et lecture :**

12 millions 457 mille 892

---

**IV. Conclusion**

A travers la pagination d'un manuscrit arabe algérien du début du 19$^{ème}$ siècle nous redécouvrons la forme originale que les chiffres arabes avaient avant de passer en Europe et de subir les transformations qui ont donné les chiffres arabes modernes. Cette forme originale, qui a complètement disparu, montre que les chiffres arabes ont pour origine les lettres arabes. Contrairement à ce que prétendent certaines hypothèses, particulièrement celles qui les présentent comme dérivant de caractères indiens, les 10 chiffres arabes que nous utilisons ne sont en fait que 10 lettres arabes plus ou moins modifiées données dans un ordre bien déterminé qui est l'ordre "Abjadi". L'hypothèse de l'origine indienne des chiffres arabes se révèle une grosse erreur démentie par la forme des chiffres arabes et par la logique de la sociologie droite-gauche dans la représentation (lecture et écriture) des nombres et dans les algorithmes des opérations élémentaires. Cette logique est toujours en vigueur chaque fois qu'il s'agit de traiter avec des nombres. Les chiffres arabes qui ont beaucoup simplifié l'écriture des nombres et les algorithmes des opérations élémentaires ont été très utilisés sous leur forme originale, celle du manuscrit, dans le Ghoubar (calcul) et dans les dates, au Maghreb (Algérie) jusqu'au début du 19$^{ème}$ siècle, ce qui prouve qu'ils sont nés au Maghreb. De Béjaïa (Bougie) ils sont passés en Europe pour donner après évolution les chiffres arabes modernes 0, 1, 2, 3, 4, 5, 6, 7, 8, 9. Ils sont aussi passés au Moyen Orient pour donner après transformations et en ajoutant deux lettres hébraïques les chiffres arabes "Mashreki" utilisés actuellement au Moyen Orient ٠ ١ ٢ ٣ ٤ ٥ ٦ ٧ ٨ ٩ .
Pour des raisons historiques, décadence et colonisation, les chiffres arabes sont devenus, durant tout le 19$^{ème}$ siècle et la première moitié du 20$^{ème}$ siècle, étrangers et inconnus dans leur propre pays, le Maghreb. Aujourd'hui les chiffres arabes modernes revenus d'Europe sont les chiffres officiels dans les pays du Maghreb (Algérie, Maroc et Tunisie).

**Références**

**Remerciements**